\documentclass[12pt]{amsart}
\usepackage{amssymb}
\usepackage{latexsym}
\usepackage{amsmath}
\usepackage{euscript}
\usepackage{graphics}

\usepackage{amssymb}
\usepackage{latexsym}
\usepackage{amsmath,amsthm}
\usepackage{graphics}
\usepackage[english]{babel}

      \def\sC{{\mathfrak C}}
      
   \def\sH{{\mathfrak H}}

\def\sS{{\mathfrak S}}

      \def\dC{{\mathbb C}}

   \def\dN{{\mathbb N}}   
      \def\dR{{\mathbb R}}

   \def\dZ{{\mathbb Z}}

   \def\cH{{\mathcal H}}

\def\cP{{\mathcal P}}      
      
   \def\cW{{\mathcal W}}

\def\wt#1{{{\widetilde #1} }}

\def\bm\chi{\mbox{\boldmath$\chi$}}

\def\ran{{\rm ran\,}}

\def\e{{\rm e}}
\let\xker=\ker \def\ker{{\xker\,}}

\def\supp{{\rm supp\,}}

\def\e{\varepsilon}
\def\wB{\widetilde{B}}

\def\f{\varphi}
\newcommand{\sddots}{\begin{picture}(2,2)
\multiput(0,0)(1.5,1){3}{.}\end{picture}}

\def\sg{\operatorname{sign}}
\def\max{\operatorname{max}}
\def\spn{\operatorname{span}}

\def\deg{\operatorname{deg}}

\def\deg{\operatorname{deg}}

\def\dist{\operatorname{dist}}

\def\wB{\widetilde{B}}

\newtheorem{theorem}{Theorem}[section]
\newtheorem{proposition}[theorem]{Proposition}

\newtheorem{lemma}[theorem]{Lemma}
\newtheorem{definition}[theorem]{Definition}
\theoremstyle{definition}
\newtheorem{example}[theorem]{Example}
\newtheorem{remark}[theorem]{Remark}

\numberwithin{equation}{section}

\numberwithin{equation}{section}

\begin{document}
\title[Generalized Jacobi operators in Krein spaces]
{Generalized Jacobi operators in Krein spaces}
\author{Maxim Derevyagin}

\address{Department of Nonlinear Analysis \\
Institute of Applied Mathematics and Mechanics\\
R.Luxemburg str. 74 \\
83114 Donetsk, Ukraine \\}

\email{derevyagin.m@gmail.com}

\date{\today}

\begin{abstract}
A special class of generalized Jacobi operators which are
self-adjoint in Krein spaces is presented. A description of the
resolvent set of such operators in terms of solutions of the
corresponding recurrence relations is given. In particular,
special attention is paid to the periodic generalized Jacobi
operators. Finally, the spectral properties of generalized Jacobi
operators are applied to prove convergence results for Pad\'e
approximants.
\end{abstract}


\maketitle

\section{Introduction}

Let $\mu$ be a positive Borel measure having infinite support
$\supp \mu \subset[a,b]\subset\dR$. To every such a measure $\mu$
there corresponds a linear functional defined on the linear space
$\cP=\spn\{\lambda^k:\,k\in\dZ_+:=\dN\cup\{0\}\}$ by the formula
\[
s_k=\sS(\lambda^k):=\int_{a}^{b}t^kd\mu(t),
\,\,k\in\dZ_+:=\dN\cup\{0\}.
\]
The functional $\sS$ is positive definite on $\cP$, that is,
$\det(s_{i+j})_{i,j=0}^{n}>0$ for all $n\in\dZ_+$ . Besides, the
measure $\mu$ (or, equivalently, the functional $\sS$) generates
the holomorphic function
\begin{equation}\label{mark_asymp}
\widehat{\mu}(\lambda)=\sS_t\left(\frac{1}{t-\lambda}\right)=\int_{a}^{b}\frac{d\mu(t)}{t-\lambda}
=-\frac{s_0}{\lambda}-\frac{s_1}{\lambda^2}-\dots-\frac{s_n}{\lambda^{n+1}}-\dots\quad
(|\lambda|>R),
\end{equation}
where $R$ is large enough. By using the Euclidean algorithm, P.L.
Tchebyshev~\cite{Tcheb} expanded the function $\widehat{\mu}$ into
the following continued fraction
\begin{equation}\label{Jfraction}
\widehat{\mu}(\lambda)\sim
-\frac{1}{\lambda-a_0-\displaystyle{\frac{b_0^2}{\lambda-a_1-\displaystyle{\frac{b_1^2}{\ddots}}}}}=
-\frac{1}{\lambda-a_0}
\begin{array}{l} \\ -\end{array}
\frac{b_0^2}{\lambda-a_1}
\begin{array}{ccc} \\ - \end{array}
\frac{b_{1}^2}{\lambda-a_2}\begin{array}{l} \\
-\cdots \end{array},
\end{equation}
where $a_j$ are real numbers, $b_j$ are positive numbers, and
$\sup\limits_{j\in\dZ_+}\{|a_j|+b_j\}<\infty$. Such continued
fractions are called $J$-fractions~\cite{JTh}. Note that the
coefficients $a_j$ and $b_j$ are uniquely determined by the
coefficients $s_j$ of the Taylor series at infinity
(see~\eqref{mark_asymp}).

It is well-known (see~\cite{Baker},~\cite{JTh},~\cite{NikSor}),
that the  $n$-th convergent $-Q_n(\lambda)/P_n(\lambda)$ of
the continued fraction~\eqref{Jfraction} is characterized by the
following property
\begin{equation}\label{eq:1.12}
-\frac{Q_n(\lambda)}{P_n(\lambda)}=-\frac{s_{0}}{\lambda}-\frac{s_{1}}{\lambda^{2}}-\dots-\frac{s_{2n-1}}{\lambda^{2n}}
+O\left(\frac{1}{\lambda^{2n+1}}\right) \quad(\lambda{\rightarrow
}\infty).
\end{equation}
In other words, the rational function $-Q_n/P_n$ is the $n$-th
diagonal Pad\'e approximant to  $\widehat{\mu}$, that is,
\[
\widehat{\mu}(\lambda)+\frac{Q_n(\lambda)}{P_n(\lambda)}=O\left(\frac{1}{\lambda^{2n+1}}\right)\quad(\lambda{\rightarrow
}\infty).
\]
Further, by using standard argumentation (see~\cite{JTh}), one can
see that the polynomials $P_{n}$ and $Q_{n}$ are solutions of the
three-term recurrence relations
\begin{equation}\label{eq:2.10_d}
b_{j-1}u_{j-1}+a_ju_{j}+b_{j}u_{j+1}=\lambda u_j\quad(j\in\dN),
\end{equation}
with initial conditions
\begin{equation}\label{DiffEq_d}
    P_0(\lambda)=1,\quad P_1(\lambda)=\frac{p_0(\lambda)}{b_0},\quad
     Q_0(\lambda)=0,\quad Q_1(\lambda)=\frac{1}{b_0}.
\end{equation}

On the other hand, to the recurrence relations~\eqref{eq:2.10_d}
(or, equivalently, to the continued fraction~\eqref{Jfraction})
there corresponds a linear bounded self-adjoint operator  in the
space $\ell^2_{[0,\infty)}$. More precisely, that operator is
generated by the following tridiagonal matrix
\[
H=\left(%
\begin{array}{cccc}
  a_0 & b_0 &  &  \\
  b_0 & a_1 & b_1 &  \\
      & b_1 & a_2 & \ddots \\
      &     & \ddots & \ddots \\
\end{array}%
\right)
\]
and its $m$-function
$m(\lambda)=\left((H-\lambda)^{-1}e,e\right)_{\ell^2_{[0,\infty)}}$,
where $e=(1,0,\dots,0,\dots)^\top\in\ell^2_{[0,\infty)}$,
coincides with $\widehat{\mu}(\lambda)$. So, one has an operator
representation of $\widehat{\mu}$
\[
\widehat{\mu}(\lambda)=\left((H-\lambda)^{-1}e,e\right)_{\ell^2_{[0,\infty)}},\quad
\lambda\in\dC\setminus\dR.
\]

The relations between orthogonal polynomials, Pad\'e approximants,
and Jacobi operators are well-known. These relations allow to use
operator methods to the investigation of orthogonal polynomials
and the Pad\'e approximants.

Note, that the above-mentioned results are also valid if the
underlying functional $\sS$ is not positive and has the property
\begin{equation}\label{reg_fun}
\det(s_{i+j})_{i,j=0}^{n}\ne 0,\quad s_j=\sS(\lambda^j)
\end{equation}
for all $n\in\dZ_+$ (see~\cite{AptKV},~\cite{Becker},~\cite{BK}).
The functional $\sS$ having the property~\eqref{reg_fun} is called
{\it regular}.

In the present paper a similar relations for non-regular
functionals are presented. In fact, the scheme proposed
in~\cite{AptKV},~\cite{BK} for investigation of the convergence of
Pad\'e approximants is generalized here to the non-regular case.
The paper is organized as follows. In Section~2, starting from a
(not necessarily regular) functional on $\cP$, three-term
recurrence relations for associated polynomials are derived. In
Section~3, a special class of generalized Jacobi matrices is
presented. Weyl solutions of the three-term recurrence relations
and Weyl functions are introduced in Section~4. In Section~5,
following the scheme proposed in~\cite{AptKV}, the
characterization of resolvent sets of generalized Jacobi operators
is obtained. Section~6 is concerned with the case of periodic
generalized Jacobi matrices. In Section~7, convergence results for
Pad\'e approximants are proved.


\section{Associated polynomials}

In this section, starting from a (not necessarily regular)
functional on $\cP$, three-term recurrence relations for
associated polynomials are derived.

Let us consider a holomorphic in a neighborhood of infinity
function $\f$ such that
\[
\f^{\sharp}(\lambda):=\overline{\f(\overline{\lambda})}=\f(\lambda).
\]
So, $\f$ has the Taylor expansion at infinity
\[
\displaystyle{\f(\lambda)=-\sum_{j=0}^{\infty}\frac{s_j}{\lambda^{j+1}}},
\]
where $s_j\in\dR$.  To every such a function one can associate a
real linear functional on $\cP$ defined by the formula
\[
\sS(\lambda^k):=\frac{1}{2\pi i}\oint_{|\lambda|=R}
\lambda^k\f(\lambda)d\lambda=s_k\in\dR,
\,\,k\in\dZ_+:=\dN\cup\{0\},
\]
for sufficiently large $R$. Clearly, the functional $\sS$ is not
necessarily regular, that is $\det(s_{i+j})_{i,j=0}^{n}=0$ may
vanish for some $n\in\dZ_+$ (for instance, see~\cite{DD}). In
general case, the functional $\sS$ generates an indefinite inner
product on $\cP$ (see~\cite{vRos},~\cite{AI})
\[
[f,g]_{\sS}:=\sS(f(\lambda)g^{\sharp}({\lambda}))=\frac{1}{2\pi
i}\oint_{|\lambda|=R}f(\lambda)g^{\sharp}({\lambda})\f(\lambda)d\lambda,
\quad f,g\in\cP,
\]
which is degenerate if and only if $\varphi$ is rational. In what
follows we suppose that $\varphi$ is not rational. As in the
regular case, one can associate to $\sS$ the following holomorphic
function
\begin{equation}\label{s_series}
\sS_z\left(\frac{1}{z-\lambda}\right)=\frac{1}{2\pi
i}\oint_{|z|=R} \frac{\f(z)dz}{z-\lambda}
=-\frac{s_0}{\lambda}-\frac{s_1}{\lambda^2}-\dots-\frac{s_n}{\lambda^{n+1}}-\dots\quad(|\lambda|>R).
\end{equation}

 Throughout this paper we suppose that the sequence ${\bf
s}:=\{s_j\}_{j=0}^\infty$ is {\it normalized}, i.e. the first
nonvanishing moment has modulus 1. A number $n_j\in\dN$ is called
{\it a normal index} if $\det (s_{i+k})_{i,k=0}^{n_j-1}\ne 0$.
Since $\varphi$ is not rational, there exists an infinite number
of normal indices (see~\cite[Section~16.10.2]{G}). Let
$n_1<n_2<\dots<n_{j}<\dots$ be a sequence of all normal indices.
By the choice of $n_1$ one has $s_{n_1-1}\ne 0$. Let us set
$\varepsilon_0=s_{n_1-1}$ ($|\varepsilon_0|=1$) and $\f_0:=\f$.
The principal part of the Laurent expansion for ${\displaystyle
-\frac{1}{\f_0}}$ is a polynomial of degree $k_0:=n_1$ with the
leading coefficient $\varepsilon_0$. So, we have
\begin{equation}\label{Q1}
   -\frac{1}{\f_0(\lambda)}=\varepsilon_0p_0(\lambda)+
   b_0^2\f_1(\lambda),
    \quad b_0>0
\end{equation}
where $p_0$ is a monic polynomial of degree $k_0$ and  $\f_1$ is
holomorphic in a neighborhood of infinity. Furthermore, the
function  $\f_1$ satisfies the relation
$\f_1^{\sharp}({\lambda})=\f_1(\lambda)$. Choose $b_0>0$ such that
the sequence ${\bf s}^{(1)}=\{s_j^{(1)}\}_{j=0}^\infty$ defined by
the following expansion
\begin{equation}\label{Q1asympt}
\f_1(\lambda)=
-\frac{s_{0}^{(1)}}{\lambda}-\frac{s_{1}^{(1)}}{\lambda^{2}}-\dots
-\frac{s_{2(n-k_0)}^{(1)}}{\lambda^{2(n-k_0)+1}}-\dots,
\end{equation}
at $\infty$ is normalized. This completes the first step of
expanding the series~\eqref{s_series} into a continued fraction
(see~\cite{Mag1}).

As was shown in~\cite{De}, the set of the normal indices of ${\bf
s}^{(1)}$ coincides with the following sequence
\[
n_2-k_0<\dots<n_{j}-k_0<\dots.
\]
Now one can apply the above reasoning to the function $\f_1$ and
so on. By recursion we obtain the following P-fraction
\begin{equation}\label{ContF2}
-\frac{\varepsilon_0}{p_0(\lambda)}
\begin{array}{l} \\ - \end{array}
\frac{\varepsilon_0\varepsilon_1b_0^2}{p_1(\lambda)}
\begin{array}{ccc} \\ - & \cdots & -\end{array}
\frac{\varepsilon_{j-1}\varepsilon_j b_{j-1}^2}{p_j(\lambda)}
\begin{array}{cc} \\ - & \cdots \end{array},
\end{equation}
where $\e_j=\pm 1$, $b_j>0$ and
$p_{j}(\lambda)={\lambda}^{k_{j}}+p_{k_{j}-1}^{(j)}{\lambda}^{k_{j}-1}+\dots+
p_{1}^{(j)}\lambda+p_{0}^{(j)}$ are real monic polynomials of
degree ${k_j}$ (see also~\cite{De},~\cite{DD1}). Note, that
$n_j=k_0+k_1+\dots+k_{j-1}$.

It also should be mentioned that there exist functions for which
the set $ \{b_j, p_{0}^{(j)}, \dots, p_{k_{j}-1}^{(j)}:
j\in\dZ_+\}$ of coefficients of the P-fraction is not necessarily
bounded. In particular, the Cauchy transform of the signed measure
constructed in~\cite{Stahl} gives such an example.

The continued fraction~\eqref{ContF2} can be considered as a
sequence of the linear-fractional transformations
(see~\cite[Section 5.2]{JTh})
\[
T_j(\omega):=\frac{-\e_j}{p_j(\lambda)+\e_jb_j^2\omega}
\]
having the following matrix representation
\begin{equation}\label{Wj}
\cW_j(\lambda)=\begin{pmatrix}0 & -\frac{\varepsilon_j}{b_j}\\
                            \varepsilon_jb_j &  \frac{p_j(\lambda)}{b_j}
                            \end{pmatrix},\quad j\in\dZ_+.
\end{equation}
The superposition $T_{0}\circ T_{1}\circ\dots \circ T_{j}$ of the
linear-fractional transformations corresponds to the product of
the matrices $\cW_l(\lambda)$
\begin{equation}\label{W}
\cW_{[0,j]}(\lambda)=(w_{ik}^{(j)}(\lambda))_{i,k=1}^2:=\cW_0(\lambda)\cW_1(\lambda)\dots
\cW_j(\lambda).
\end{equation}

To give an explicit formula for $\cW_{[0,j]}$ in terms of  $p_j$,
$b_j$, $\varepsilon_j$, define the polynomials
$P_{j+1}(\lambda)$, $Q_{j+1}(\lambda)$ by the equalities
\begin{equation}\label{PQj}
   \left(%
\begin{array}{c}
  -Q_0 \\
  P_0 \\
\end{array}%
\right)=\left(%
\begin{array}{c}
  0 \\
  1 \\
\end{array}%
\right),\quad
\left(%
\begin{array}{c}
 -Q_{j+1}(\lambda) \\
 P_{j+1}(\lambda) \\
\end{array}%
\right):=\cW_{[0,j]}(\lambda)\left(%
\begin{array}{c}
  0 \\
  1 \\
\end{array}%
\right), \quad j\in\dZ_+ .
\end{equation}
The relation
$\cW_{[0,j]}(\lambda)=\cW_{[0,j-1]}(\lambda)\cW_{j}(\lambda)$
(see~\eqref{W}) yields
\begin{equation}\label{W0j}
\cW_{[0,j]}(\lambda)\left(%
\begin{array}{c}
  1 \\
  0 \\
\end{array}%
\right)=
\cW_{[0,j-1]}(\lambda)\left(%
\begin{array}{c}
  0 \\
  \varepsilon_jb_j \\
\end{array}%
\right)=
\left(%
\begin{array}{c}
  -\varepsilon_jb_j Q_j(\lambda) \\
  \varepsilon_jb_j P_j(\lambda) \\
\end{array}%
\right), \quad j\in\dN .
\end{equation}
So, the matrix $\cW_{[0,j]}(\lambda)$ has the form
\begin{equation}\label{SolMatr}
\cW_{[0,j]}(\lambda)=\left(%
\begin{array}{cc}
  -\varepsilon_jb_j Q_j(\lambda) & -Q_{j+1}(\lambda) \\
  \varepsilon_jb_j P_j(\lambda) & P_{j+1}(\lambda) \\
\end{array}%
\right), \quad j\in\dZ_+ .
\end{equation}
Further, the equality
\[
\left(%
\begin{array}{c}
 -Q_{j+1}(\lambda) \\
 P_{j+1}(\lambda) \\
\end{array}%
\right)=\cW_{[0,j-1]}(\lambda)\cW_{j}(\lambda)\left(%
\begin{array}{c}
  0 \\
  1 \\
\end{array}%
\right)=
\frac{1}{b_j}\cW_{[0,j-1]}(\lambda)\left(%
\begin{array}{c}
  -\varepsilon_j \\
  p_j(\lambda) \\
\end{array}%
\right), \quad j\in\dN,
\]
 shows that the polynomials $P_{j}(\lambda)$, $Q_{j}(\lambda)$ are solutions of the difference equation
\begin{equation}\label{eq:2.10}
\varepsilon_{j-1}\varepsilon_j
b_{j-1}u_{j-1}-p_j(\lambda)u_{j}+b_{j}u_{j+1}=0\,\,\,(j\in\dN),
\end{equation}
obeying the initial conditions
\begin{equation}\label{DiffEq}
 \begin{split}
    P_0(\lambda)&=1,\quad P_1(\lambda)=\frac{p_0(\lambda)}{b_0},\\
     Q_0(\lambda)&=0,\quad Q_1(\lambda)=\frac{\varepsilon_0}{b_0}.
\end{split}
\end{equation}
According to~\eqref{SolMatr}, the $(j+1)$-th convergent of the
continued fraction~\eqref{ContF2} is equal to
\[
f_j:=T_{0}\circ T_{1}\circ\dots \circ
T_{j}(0)=-Q_{j+1}(\lambda)/P_{j+1}(\lambda).
\]
The relations~\eqref{SolMatr}, \eqref{W}, and \eqref{Wj} imply the
following statement.
\begin{proposition}[\cite{DD1}]\label{Ostrogr}
The polynomials $P_j$, $Q_j$ satisfy the following generalized
Liouville-Ostrogradsky formula
\begin{equation}\label{Ostrogr2}
\varepsilon_jb_j(Q_{j+1}\left(\lambda)P_{j}(\lambda) -
        Q_{j}(\lambda) P_{j+1}(\lambda) \right)=1\quad(j\in\dZ_+).
\end{equation}
\end{proposition}

\section{Generalized Jacobi matrices}

The main goal of this section is to present a special class of
generalized Jacobi matrices.

 Let $p(\lambda)=p_{n}{\lambda}^{n}+\dots+p_{1}\lambda+p_{0}$
be a monic scalar real polynomial of degree $n$, i.e. $p_{n}=1$.
Let us associate to the polynomial $p$ its symmetrizator  $E_p$
and let the companion matrix $C_p$ be given by
\begin{equation}\label{comp}
E_{p}=\begin{pmatrix}
p_{1}&\dots&p_{n}\\
\vdots&\sddots&\\
p_{n}&&{\bf 0}\\
\end{pmatrix},\quad
C_{p}=\begin{pmatrix}
0&\dots&0&-p_{0}\\
1&&{\bf 0}&-p_{1}\\
&\ddots&&\vdots\\
{\bf 0}&&1&-p_{n-1}\\
\end{pmatrix}.
\end{equation}
As is known, $\det(\lambda-C_p)=p(\lambda)$ and the spectrum
$\sigma(C_p)$ of the companion matrix $C_p$ is simple. The
matrices $E_{p}$ and $C_{p}$ are related by (see~\cite{GLR})
\begin{equation}\label{cb}
C_{p}E_{p}=E_{p}C_{p}^{\top}.
\end{equation}
So, $C_{p}E_{p}$ is a symmetric matrix.

\begin{definition}[\cite{DD}, \cite{KL79}]
Let $p_j$ be real monic polynomials of degree ${k_j}$
\[
p_{j}(\lambda)={\lambda}^{k_{j}}+p_{k_{j}-1}^{(j)}{\lambda}^{k_{j}-1}+\dots+
p_{1}^{(j)}\lambda+p_{0}^{(j)},
\]
and let  $\e_{j}=\pm 1$,  $b_{j}>0$ $(j\in\dN)$. The tridiagonal
block matrix
\begin{equation}\label{mJacobi}
H=\begin{pmatrix}
A_{0}   &\wB_{0}&       &{\bf 0}\\
B_{0}   &A_1    &\wB_{1}&\\
        &B_1    &A_{2} &\ddots\\
{\bf 0} &       &\ddots &\ddots\\
\end{pmatrix}
\end{equation}
where $A_{j}=C_{p_{j}}$ and $k_{j+1}\times k_{j}$ matrices $B_{j}$
and $k_{j}\times k_{j+1}$ matrices $\wB_{j}$ are given by
\begin{equation}\label{bblock}
B_{j}=\begin{pmatrix}
0&\dots&b_{j}\\
\hdotsfor{3}\\
0&\dots&0\\
\end{pmatrix},\,
\wB_{j}= \begin{pmatrix}
0&\dots&\wt b_{j}\\
\hdotsfor{3}\\
0&\dots&0\\
\end{pmatrix} \,
(\wt b_{j}=\e_j\e_{j+1} b_j, \, j=0,\dots, N-1),
\end{equation}
will be called a {\it generalized Jacobi matrix}  associated with
the sequences of polynomials $\{{\e_{j}p_{j}\}}_{j=0}^{\infty}$
and numbers $\{b_j\}_{j=0}^{\infty}$.
\end{definition}
\begin{remark}
The papers~\cite{DD}, \cite{DD1}, and~\cite{KL79} are only
concerned with the case of generalized Jacobi matrices which are
finite rank perturbations of classical Jacobi matrices. In fact,
the generalized Jacobi matrix in question is associated to the
P-fraction~\eqref{ContF2} or, equivalently, to the sequence of
matrices $\cW_j$ having the form~\eqref{Wj}.
\end{remark}

\renewcommand{\theequation}{A.\arabic{equation}}
\setcounter{equation}{0}

{\bf From now on, we suppose that}
\begin{equation}\label{assumption1}
\text{there exists }N\in\dN:\quad\deg p_j\le N,\quad j\in\dZ_+,
\end{equation}
\begin{equation}\label{assumption2}
\sup \{b_j, |p_{k}^{(j)}|: j\in\dZ_+, k=0,\dots,k_j-1\}< +\infty.
\end{equation}
\renewcommand{\theequation}{\thesection.\arabic{equation}}
\setcounter{equation}{4}

Let $\ell^2_{[0,\infty)}$ denote the Hilbert space of complex
square summable sequences $(w_0,w_1,\dots)$ with the usual inner
product. Setting
\begin{equation}\label{NormInd}
n_{0}=0, \quad n_j=\sum_{i=0}^{j-1}k_{i}\quad(j\in\dN),
\end{equation}
define a standard basis in $\ell^2_{[0, \infty)}$ by the
equalities
\[
e_{j,k}=\{{\delta_{l, n_{j}+k}\}}_{l=0}^{\infty} \quad (j\in\dZ_+;
k=0,\dots,k_{j}-1),\quad e:=e_{0,0}.
\]
 Define the symmetric matrix
$G$ by the equality
\begin{equation}\label{Gram}
G=\mbox{diag}(G_{0},G_1,\dots),\quad
G_{j}=\e_{j}E_{p_{j}}^{-1}\quad (j\in\dZ_+).
\end{equation}
Further, we may identify via usual matrix product the matrix $G$
with an operator on the linear space $\sC_0$ of finite sequences
of $\ell^2_{[0, \infty)}$. Its closure will be also denoted by
$G$. In view of~\eqref{assumption1},~\eqref{assumption2}, the
operator $G$ defined on $\ell^2_{[0,\infty)}$ is bounded and
self-adjoint. Moreover, $G^{-1}$ is a bounded linear operator in
$\ell^2_{[0,\infty)}$.

 Let $H_{[j,l]}$ ($G_{[j,l]}$) be a submatrix of $H$ ($G$), corresponding to the
basis vectors $\{e_{i,k}\}^{k=0,\dots,k_{i}-1}_{i=j,\dots,l}$
$(0\le j\le l < +\infty)$. The matrix $H_{[j,l]}$ will be called a
finite generalized Jacobi matrix.

 Let $\sH_{[0,\infty)}$ be a space of elements of $\ell^2_{[0,\infty)}$
 provided with the following indefinite inner product
\begin{equation}\label{metricInfty}
\left[x,y\right]=(Gx,y)_{\ell^2_{[0,\infty)}}\quad
(x,y\in\ell^2_{[0,\infty)}).
\end{equation}

Let us recall~\cite{AI} that a pair $(\sH, [\cdot,\cdot])$
consisting of a Hilbert space $\sH$ and a sesquilinear form
$[\cdot,\cdot]$ on $\sH\times\sH$ is called {\it a space with
indefinite inner product}. A space with indefinite metric $(\sH,
[\cdot,\cdot])$ is called {\it a Krein space} if the indefinite
scalar product $[\cdot,\cdot]$ can be represented as follows
\[
[x,y]=(Jx,y)_{\sH}\quad x,y\in\sH,
\]
where the linear operator $J$ satisfies the following conditions
\[
J=J^{-1}=J^*.
\]
The operator $J$ is called {\it the fundamental symmetry}. So, one
can see that the space $\sH_{[0,\infty)}$ is the Krein space with
the fundamental symmetry $J=\sg G$ (see~\cite{AI} for details).

\begin{proposition}\label{SymGJM}
Under the assumptions~\eqref{assumption1},~\eqref{assumption2},
the considered generalized Jacobi matrix defines a bounded
self-adjoint operator $H$ (a generalized Jacobi operator) in the
Krein space $\sH_{[0,\infty)}$, that is,
\begin{equation}\label{hermitian}
\left[H x,y\right]=\left[x,H y\right]\quad x,y\in\sH_{[0,\infty)}.
\end{equation}
\end{proposition}

\begin{proof} It is not hard to see that, according to~\eqref{assumption1},~\eqref{assumption2},
the matrix in question generates a bounded operator in
$\sH_{[0,\infty)}$. Relation~\eqref{hermitian} is implied by
~\eqref{cb} (see~\cite{DD} for details).
\end{proof}

  Let us extend the system
$\{P_j(\lambda)\}_{j=0}^{\infty}$,
$\{Q_j(\lambda)\}_{j=0}^{\infty}$ by the equalities
\begin{equation}\label{polynom0}
P_{j,k}(\lambda)={\lambda}^{k}P_{j}(\lambda),\quad
Q_{j,k}(\lambda)={\lambda}^{k}Q_{j}(\lambda)\quad
(j\in\dZ_+;k=0,\dots,k_{j}-1).
\end{equation}
Setting
\[
{\bf
P}_{[l,j]}(\lambda)=(P_{l,0}(\lambda),\dots,P_{l,k_l-1}(\lambda),\dots,
P_{j,0}(\lambda),\dots,P_{j,k_j-1}(\lambda)),
\]
\[
{\bf
Q}_{[l,j]}(\lambda)=(Q_{l,0}(\lambda),\dots,Q_{l,k_l-1}(\lambda),\dots,
Q_{j,0}(\lambda),\dots,Q_{j,k_j-1}(\lambda)).
\]
one can rewrite the system~\eqref{eq:2.10}, \eqref{polynom0} in
the following manner
\begin{equation}\label{Vector10}
{\bf
P}_{[0,j]}(\lambda)(\lambda-H_{[0,j]})=(0,\dots,0,b_jP_{j+1,0}(\lambda))\quad
(j\in\dZ_+),
\end{equation}
\begin{equation}\label{Vector11}
{\bf
Q}_{[0,j]}(\lambda)(\lambda-H_{[0,j]})=(\underbrace{0,\dots,0,-\varepsilon_0}_{k_0},0,\dots,0,b_jQ_{j+1,0}(\lambda))\quad
(j\in\dZ_+).
\end{equation}
Since $Q_{0,0}(\lambda)=\dots=Q_{0,k_0-1}(\lambda)\equiv 0$, the
relation~\eqref{Vector11} reduces to
\begin{equation}\label{Vector12}
{\bf
Q}_{[1,j]}(\lambda)(\lambda-H_{[1,j]})=(0,\dots,0,b_jQ_{j+1,0}(\lambda))\quad
(j\in\dN),
\end{equation}
It follows from~\eqref{Vector10} and \eqref{Vector12} that the
eigenvalues of $H_{[0,j]}$ and $H_{[1,j]}$ coincide with the roots
of $P_{j+1}(\lambda)$ and $Q_{j+1}(\lambda)$, respectively.

\begin{proposition}[\cite{DD}]\label{detf}
The polynomials $P_{j}$ and $Q_j$  $(j\in\dN)$ can be found by the
formulas
\begin{eqnarray}
P_{j}(\lambda)&=&(b_{0}\dots b_{j-1})^{-1}\det(\lambda-H_{[0,j-1]}),\label{polynom1}\\
Q_{j}(\lambda)&=&\varepsilon_0(b_{0}\dots
b_{j-1})^{-1}\det(\lambda-H_{[1,j-1]}).\label{polynom2}
\end{eqnarray}
\end{proposition}
The formulas~\eqref{polynom1} and~\eqref{polynom2} in the
classical case can be found in~\cite[Section 7.1.2]{Be}
and~\cite[Section 6.1]{Atk}. The following statement is an easy
consequence of the recurrence relations~\eqref{eq:2.10}.
\begin{proposition}[\cite{DD},\cite{DD1}]\label{prostota}
Let $j\in\dN$. Then
\begin{enumerate}
\item[i)] The polynomials $P_j$ and $P_{j+1}$  have no common zeros.
\item[ii)] The polynomials $Q_j$ and $Q_{j+1}$ have no common zeros.
\item[iii)] The polynomials $P_j$ and $Q_j$ have no common zeros.
\end{enumerate}
\end{proposition}

Taking into account the equality
$G_{[0,j]}H_{[0,j]}=H_{[0,j]}^\top G_{[0,j]}$ which is implied
by~\eqref{cb} (see~\cite{DD}) and setting
\[
\pi_{[0,j]}(\lambda)=G_{[0,j]}^{-1}{\bf
P}_{[0,j]}(\lambda)^\top,\quad
\xi_{[0,j]}(\lambda)=G_{[0,j]}^{-1}{\bf Q}_{[0,j]}(\lambda)^\top,
\]
one can rewrite~\eqref{Vector10}, \eqref{Vector11} in the form
\begin{equation}\label{Matrix10}
(\lambda-H_{[0,j]})\pi_{[0,j]}(\lambda)=
\varepsilon_jb_jP_{j+1,0}(\lambda)e_{j,0}\quad (j\in\dZ_+),
\end{equation}
\begin{equation}\label{Matrix11}
(\lambda-H_{[0,j]})\xi_{[0,j]}(\lambda)+e_{0,0}=
\varepsilon_jb_jQ_{j+1,0}(\lambda)e_{j,0}\quad (j\in\dZ_+).
\end{equation}
Further, let us set
\begin{equation}\label{defpiksi}
\pi(\lambda)=G^{-1}(P_{0,0}(\lambda),\dots,P_{0,{k_0}-1}(\lambda),\dots)^{\top},\quad
\xi(\lambda)=G^{-1}(Q_{0,0}(\lambda),\dots,Q_{0,{k_0}-1}(\lambda),\dots)^{\top}.
\end{equation}
Now, it follows from~\eqref{Matrix10}-\eqref{defpiksi} that the
following formal equalities hold true
\begin{equation}\label{formal_eq}
(\lambda-H)\pi(\lambda)=0,\quad (\lambda-H)\xi(\lambda)=-e_{0,0}.
\end{equation}

The first equality in~\eqref{formal_eq} allows us to characterize
the point spectrum of  $H$.

\begin{proposition}
$\lambda\in\sigma_p(H)$ if and only if
$\pi(\lambda)\in\ell^2_{[0,\infty)}$.
\end{proposition}
By the definition of $\pi(\lambda)$ and the
assumptions~\eqref{assumption1},~\eqref{assumption2}, we see that
\[
\pi(\lambda)\in\ell^2_{[0,\infty)}\Longleftrightarrow
\sum_{j=0}^{\infty}|P_{j}(\lambda)|^2<+\infty.
\]

\section{Weyl solutions and Weyl functions}

If for some $\lambda\in\dC$ there exists a solution
$\{W_j(\lambda)\}_{j=0}^{\infty}$ of the recurrence
relations~\eqref{eq:2.10} such that
\begin{equation}\label{Weylsolution}
\{W_j(\lambda)\}_{j=0}^{\infty}\in\ell^2_{[0,\infty)}\text{ and
}\{W_j(\lambda)\}_{j=0}^{\infty}\ne\{P_j(\lambda)\}_{j=0}^{\infty}
\end{equation}
then we will say that there exists a Weyl solution
$\{W_j(\lambda)\}_{j=0}^{\infty}$ of the recurrence
relations~\eqref{eq:2.10} at the point $\lambda$. Since
$\{P_j(\lambda)\}_{j=0}^{\infty}$ and
$\{Q_j(\lambda)\}_{j=0}^{\infty}$ are linearly independent
solutions of~\eqref{eq:2.10}, the Weyl solution admits the
following representation
\begin{equation}\label{reprWs}
W_j(\lambda)=Q_{j}(\lambda)+m(\lambda)P_{j}(\lambda),
\end{equation}
where $m(\lambda)$ is a complex number. The following statement
shows the relation between $m(\lambda)$ and the operator  $H$.
\begin{proposition}\label{Wsexplicitform}
Let $\lambda\in\rho(H)$ and let
\begin{equation}\label{Wfexf}
m(\lambda)=[(H-\lambda)^{-1}e,e],\quad e:=e_{0,0}.
\end{equation}
Then the family $\{W_j(\lambda)\}_{j=0}^{\infty}$ given
by~\eqref{reprWs} is the Weyl solution of ~\eqref{eq:2.10} at
$\lambda$. Moreover, $\{W_j(\lambda)\}_{j=0}^{\infty}$ has the
form
\begin{equation}\label{Wsexf}
W_j(\lambda)=[(H-\lambda)^{-1}e,e_{j,0}],\quad j\in\dZ_+.
\end{equation}
\begin{proof}
For $\lambda\in\rho(H)$ the relation~\eqref{formal_eq} implies
that there exists a number $m(\lambda)\in\dC$ such that
\[
\xi(\lambda)+m(\lambda)\pi(\lambda)=(H-\lambda)^{-1}e\in\ell^2_{[0,\infty)}.
\]
Since $G$ is a bounded operator, we have
\begin{equation}\label{l2belong}
G(\xi(\lambda)+m(\lambda)\pi(\lambda))\in\ell^2_{[0,\infty)}.
\end{equation}
Using~\eqref{defpiksi} yields
\[
[(H-\lambda)^{-1}e,e_{j,0}]=(G(\xi(\lambda)+m(\lambda)\pi(\lambda)),e_{j,0})_{\ell^2}=Q_{j}(\lambda)+m(\lambda)P_{j}(\lambda)=:W_j(\lambda)\,(\ne
P_j(\lambda)).
\]
The latter relation means that $\{W_j(\lambda)\}_{j=0}^{\infty}$
is a solution of~\eqref{eq:2.10}. Due to~\eqref{l2belong}, we
obtain that
$\{W_j(\lambda)\}_{j=0}^{\infty}\in\ell^2_{[0,\infty)}$. Now, it
follows from~\eqref{DiffEq} that
\[
W_0(\lambda)=m(\lambda)=[(H-\lambda)^{-1}e,e].
\]
\end{proof}
\end{proposition}

\begin{definition}\label{Weylfunction}
The function $m$ defined by~\eqref{reprWs} is called a Weyl
function of the operator $H$.
\end{definition}
\begin{remark}
It should be mentioned that a general treatment of the Weyl functions
of classical Jacobi matrices in the framework of extension theory of non-densely
defined symmetric operators was proposed in~\cite{MMal92}. A general
treatment of the Weyl functions of symmetric operators in Krein spaces
in the framework of extension theory was presented in~\cite{VDer99}.
The function $m$ defined by~\eqref{Wfexf} on $\rho(H)$ is also
called the $m$-function of $H$ (see~\cite{DD}, \cite{GS}).
\end{remark}

Since $H$ is bounded,  $m$ admits the representation
\begin{equation}\label{mreprinft}
 m(\lambda)= [(H-\lambda)^{-1}e,e]=-\sum_{j=0}^{\infty}\frac{
{s_{j}}}{\lambda^{j+1}},\quad (|\lambda|>\Vert H\Vert)
\end{equation}
where ${s_{i}}=[H^ie,e]$.

Analogously, one can define $m$-functions of shortened generalized
Jacobi matrices.
\begin{definition}\label{trWf}
The function
\begin{equation}\label{SWfun}
    m_{[0,j]}(\lambda)=[(H_{[0,j]}-\lambda)^{-1}e,e]
\end{equation}
is called the $m$-function of $H_{[0,j]}$.
\end{definition}

Making use of the structure of $H_{[0,j]}$, we obtain that
(see~\cite{DD})
\begin{equation}\label{SWfun2}
    m_{[0,j]}(\lambda)=-\e_0\frac{\det(\lambda-H_{[1,j]})}{\det(\lambda-H_{[0,j]})}.
\end{equation}
According to~\eqref{polynom1}, \eqref{polynom2}, the
formula~\eqref{SWfun2} can be rewritten as follows
\begin{equation}~\label{mQP}
m_{[0,j]}(\lambda)=-\frac{Q_{j+1}(\lambda)}{P_{j+1}(\lambda)}.
\end{equation}
It follows from~\eqref{eq:2.10} (see~\cite{DD} for details) that
the $m$-function $m_{[0,j]}(\lambda)$ and the $m$-function
$m_{[1,j]}(\lambda)$ of $H_{[1,j]}$ are related by the equality
\begin{equation}\label{Riccati}
m_{[0,j]}(\lambda)=\frac{-\varepsilon_0}
    {p_0(\lambda)+\varepsilon_0b_0^2m_{[1,j]}(\lambda)},\quad j\in\dN.
\end{equation}

An analogous statement for infinite generalized Jacobi matrices is
an essential ingredient in the proof of the following result.

\begin{theorem}
Under the assumptions~\eqref{assumption1},~\eqref{assumption2},
the generalized Jacobi matrix $H$ is uniquely determined by its
Weyl function $m$.
\end{theorem}
\begin{proof}
By using the Frobenius formula, one can see that the Weyl function
$m$ of $H$ and the Weyl function
\[
m_{[1,\infty)}(\lambda)=[(H_{[1,\infty)}-\lambda)^{-1}e_{1,0},e_{1,0}]
\]
of $H_{[1,\infty)}$ are related by the equality
\begin{equation}\label{Riccatiinfty}
m(\lambda)=\frac{-\varepsilon_0}
    {p_0(\lambda)+\varepsilon_0b_0^2m_{[1,\infty)}(\lambda)},\quad |\lambda|>
\Vert H\Vert\ge\Vert H_{[1,\infty)}\Vert
\end{equation}
(a more detailed reasoning can be found in~\cite{DD}). Further,
consecutive applications of the relation~\eqref{Riccatiinfty}
leads to the $P$-fraction~\eqref{ContF2}. So, one can uniquely
recover the generalized Jacobi matrix $H$.
\end{proof}

\begin{remark}
In the definite case, formulas~\eqref{SWfun2} and~\eqref{mQP} are easy consequences 
of the theory developed in~\cite{MMal92}.
In this case, it is well known that the Weyl function determines the classical Jacobi matrix uniquely 
(for instance, see~\cite{GS}, \cite{MMal92}). 
It is worth to mention that this result as well as Borg
type uniqueness result was recently extended
to the case of normal matrices  (see~\cite{SMal}).
Besides,  a canonical form of a normal matrix
(an analog of Jacobi matrix for selfadjoint matrices)
was also introduced there.
Some inverse problems for finite generalized Jacobi matrices were
considered in~\cite{De06},~\cite{DD}.
\end{remark}

\section{The resolvent set of $H$}

Here, following the scheme proposed in~\cite{AptKV}, the
characterization of resolvent sets $\rho(H)$ of generalized Jacobi
operators is obtained. We begin with an auxiliary lemma which
gives a criterion of the density of $\ran(H-\lambda)$ in
$\ell^2_{[0,\infty)}$.

\begin{lemma}\label{lem_2_1} For $\lambda\in\dC$ the equation
\[
(H-\lambda I)x=e_{j,k}
\]
has a solution $x=x(j,k)\in\ell^2_{[0,\infty)}$ for all
$j,k\in\dZ_+$ iff there exists a Weyl solution at $\lambda$.
\end{lemma}

\begin{proof}
\noindent 1) First, let us consider the case where $j=k=0$. It
follows from~\eqref{formal_eq} that the equation
$(H-\lambda)x=e_{0,0}$ has a solution belonging to
$\ell^2_{[0,\infty)}$ if and only if there exists a number
$m(\lambda)\in\dC$  such that
$\xi(\lambda)+m(\lambda)\pi(\lambda)\in\ell^2_{[0,\infty)}$.
Furthermore, in this case we have
\begin{equation}\label{fres_e}
x=x(0,0)=\xi(\lambda)+m(\lambda)\pi(\lambda).
\end{equation}
\noindent 2) Next, let $k=0$ and $j\in\dN$. According
to~\eqref{Matrix10} and~\eqref{Matrix11}, we see that
\begin{equation}\label{f_2_2}
(H-\lambda)\xi_{[0,j]}(\lambda)=e_{0,0}-\e_jb_jQ_{j+1,0}(\lambda)e_{j,0}+\e_jb_jQ_{j,0}(\lambda)e_{j+1,0},
\end{equation}
\begin{equation}\label{f_2_3}
(H-\lambda)\pi_{[0,j]}(\lambda)=-\e_jb_jP_{j+1,0}(\lambda)e_{j,0}+\e_jb_jP_{j,0}(\lambda)e_{j+1,0},
\end{equation}
Adding~\eqref{f_2_2} multiplied by $-P_j$ and~\eqref{f_2_3}
multiplied by $Q_j$, one obtains
\begin{equation}\label{f_2_4}
(H-\lambda)\left[-P_j(\lambda)\xi_{[0,j]}(\lambda)+Q_j(\lambda)\pi_{[0,j]}(\lambda)\right]=-P_j(\lambda)e_{0,0}-\e_jb_j(P_{j+1}(\lambda)Q_j(\lambda)-Q_{j+1}(\lambda)P_j(\lambda))e_{j,0}.
\end{equation}
Due to~\eqref{Ostrogr2}, \eqref{f_2_4} can be rewritten as follows
\[
(H-\lambda)\left[-P_j(\lambda)\xi_{[0,j]}(\lambda)+Q_j(\lambda)\pi_{[0,j]}(\lambda)\right]+P_j(\lambda)e_{0,0}=e_{j,0}.
\]
The latter relation shows that the equation $(H-\lambda)x=e_{j,0}$
has a solution $x=x(j,0)\in\ell_{[0,\infty)}^{2}$ if and only if
$e_{0,0}\in\ran(H-\lambda)$. Moreover, the solution admits the
representation
\begin{equation}\label{f_2_5}
x=x(j,0)=-P_j(\lambda)\xi_{[0,j]}(\lambda)+Q_j(\lambda)\pi_{[0,j]}(\lambda)+P_j(\lambda)\left[\xi(\lambda)+m(\lambda)\pi(\lambda)\right].
\end{equation}
\noindent 3) Finally, assume that $k\ne 0$. Observe that
\begin{equation}\label{f_2_6}
(H-\lambda)e_{j,0}=e_{j,1}-\lambda e_{j,0},\dots,
(H-\lambda)e_{j,k_j-2}=e_{j,k_j-1}-\lambda e_{j,k_j-2}.
\end{equation}
The chain of equalities~\eqref{f_2_6} imply that the equation
$(H-\lambda)x=e_{j,k}$ has a solution
$x=x(j,k)\in\ell^{2}_{[0,\infty)}$ if and only if
$e_{j,0}\in\ran(H-\lambda)$ (or, equivalently,
$e_{0,0}\in\ran(H-\lambda)$). Besides, the solution $x=x(j,k)$ can
be expressed in the following manner
\begin{equation*}
x=x(j,k)=e_{j,k-1}+\lambda e_{j,k-2}+\dots+
\lambda^k(-P_j(\lambda)\xi_{[0,j]}(\lambda)+Q_j(\lambda)\pi_{[0,j]}(\lambda)+P_j(\lambda)\left[\xi(\lambda)+m(\lambda)\pi(\lambda)\right]).
\end{equation*}
\end{proof}
\begin{remark}
In fact, Lemma~\ref{lem_2_1} gives a way to express the formal
inverse operator $R_{\lambda}$ to $H-\lambda$. Moreover, it is not
so hard to see that the Weyl solution at $\lambda$ exists if and
only if $e\in\ran(H-\lambda)$.
\end{remark}

Now we are ready to prove the main result of the present paper.
\begin{theorem}\label{main_result}
Under the assumptions~\eqref{assumption1},~\eqref{assumption2},
$\lambda\in\rho(H)$ if and only if there exist a Weyl solution
$\{W_j(\lambda)\}_{j=0}^{\infty}$ at $\lambda$ and numbers
$q\in(0,1)$, $C>0$ such that
\begin{equation}\label{estimate1}
|P_i(\lambda)W_j(\lambda)|\le Cq^{n_j-n_i},\quad i\le j.
\end{equation}
\end{theorem}

\begin{proof}
Let us prove the sufficiency. Let
$\cH_k:=\overline{\spn}\{e_{j,k}|j\in\dZ_+\}$ and let
$\cP_{\cH_k}$ be the orthogonal projector onto $\cH_k$ in
$\ell^2_{[0,\infty)}$. We start with proving the boundedness  of
the operator $R_{\lambda}\cP_{\cH_k}$ for any
$k\in\{0,\dots,N-1\}$. First, it is convenient to consider the
operator $GR_{\lambda}\cP_{\cH_0}$. Taking into
account~\eqref{f_2_5}, we obtain
\[
(GR_{\lambda}\cP_{\cH_0}e_{j,0})_{i,k}=(GR_{\lambda}e_{j,0})_{i,k}:=(GR_{\lambda}e_{j,0},e_{i,k})=\begin{cases}
       \lambda^k P_i(\lambda)(Q_j(\lambda)+m(\lambda)P_j(\lambda)),& i\le j;\\
       \lambda^k P_j(\lambda)(Q_i(\lambda)+m(\lambda)P_i(\lambda)),& i>j.\\
 \end{cases}
\]
It is clear that one can represent the operator
$GR_{\lambda}\cP_{\cH_0}$ as the sum of upper and lower triangular
operators:
$GR_{\lambda}\cP_{\cH_0}=R_{\lambda}^{(1)}+R_{\lambda}^{(2)}$. To
be more precise, we choose $R_{\lambda}^{(1)}$ in the following
way
\[
R_{\lambda}^{(1)}e_{j,k}=0\text{ for }k\ne 0,\quad
R_{\lambda}^{(1)}e_{j,0}=y_{j}^{(1)},\text{ where }
\]
\[
(y_{j}^{(1)})_{i,k}=\begin{cases}
       \lambda^k P_i(\lambda)(Q_j(\lambda)+m(\lambda)P_j(\lambda)),& i\le j;\\
       0,& i>j.\\
\end{cases}
\]
Now,  one can prove that $R_{\lambda}^{(1)}$ is bounded. Indeed,
setting
\[
x=\sum\limits_{j=0}^{l}\sum\limits_{k=0}^{k_j-1}x_{j,k}e_{j,k}\text{
and }
y=\sum\limits_{j=0}^{l}\sum\limits_{k=0}^{k_j-1}y_{j,k}e_{j,k}
\]
we get the following relation
\[
\begin{split}
(R_{\lambda}^{(1)}x,y)_{\ell^2}&=\sum_{j=0}^{l}x_{j,0}\left(R_{\lambda}^{(1)}e_{j,0},\sum\limits_{i=0}^{l}\sum_{k=0}^{k_i-1}y_{i,k}e_{i,k}\right)\\&=
\sum_{j=0}^lx_{j,0}\sum_{i=0}^j\sum_{k=0}^{k_i-1}\lambda^k
P_i(\lambda)(Q_j(\lambda)+m(\lambda)P_j(\lambda))\overline{y}_{i,k}.
\end{split}
\]
Thus, \eqref{estimate1} yields
\[
|(R_{\lambda}^{(1)}x,y)_{\ell^2}|\le
C\max\{1,|\lambda|^{N-1}\}\sum_{j=0}^l|x_{j,0}|\sum_{i=0}^j\sum_{k=0}^{k_i-1}q^{n_j-n_i}|{y}_{i,k}|.
\]
Notice that $n_j-n_i=\sum\limits_{l=i+1}^{j}k_l\ge
\sum\limits_{l=i+1}^{j}1=j-i$ and, therefore, we have
\[
\begin{split}
|(R_{\lambda}^{(1)}x,y)_{\ell^2}|&\le
C\max\{1,|\lambda|^{N-1}\}\sum_{j=0}^l\sum_{i=0}^{j}\sum_{k=0}^{k_i-1}q^{{j-i}}|x_{j,0}||{y}_{i,k}|\\&\le
C\max\{1,|\lambda|^{N-1}\}\sum_{s=0}^{l}q^{{s}}\sum_{j=s}^l|x_{j,0}|\sum_{k=0}^{k_{j-s}-1}|{y}_{j-s,k}|\\
&\le C\sqrt{N}\max\{1,|\lambda|^{N-1}\}\sum_{s=0}^{l}q^{{i}}\Vert
x\Vert_{\ell^2}\Vert y\Vert_{\ell^2}
=\widetilde{C}\frac{1-q^{{l+1}}}{1-q}\Vert x\Vert_{\ell^2}\Vert
y\Vert_{\ell^2}.
\end{split}
\]
Hence, $R_{\lambda}^{(1)}$ is a bounded operator. Similarly, one
can prove the boundedness of $R_{\lambda}^{(2)}$. So, we have
proved that $GR_{\lambda}\cP_{\cH_0}$ is bounded. Since $G^{-1}$
is bounded, $R_{\lambda}\cP_{\cH_0}$ is also bounded.

Further, let $k\in\{1,\dots,N-1\}$. From~\eqref{f_2_6} one can
deduce
\[
R_{\lambda}e_{j,k}=e_{j,k-1}+\lambda
R_{\lambda}e_{j,k-1}=V_ke_{j,k}+\lambda R_{\lambda}V_ke_{j,k},
\]
where $V_k: e_{j,k}\mapsto e_{j,k-1}$ is an isometric operator
from $\cH_k$ to $\cH_{k-1}$. If $h\in\cH_k$ then $R_{\lambda}h=V_k
h+\lambda R_{\lambda}V_k h$. So, the boundedness of $V_1$ and
$R_{\lambda}\cP_{\cH_0}$ implies that $R_{\lambda}\cP_{\cH_1}$ is
bounded. Analogously, $R_{\lambda}\cP_{\cH_k}$ is bounded for
$k\in\{2,\dots,N-1\}$. This implies that
$R_{\lambda}=\sum\limits_{i=0}^{N-1}R_{\lambda}\cP_{\cH_i}$ is a
bounded operator. The latter means that the domain of
$R_{\lambda}$ is $\ell^2_{[0,\infty)}$. Since $\lambda$ is not an
eigenvalue, we have $\ker(H-\lambda)=\{0\}$ and
$\ran(H-\lambda)=\ell^2_{[0,\infty)}$. Now, applying the Banach
theorem on inverse operators we obtain that $\lambda\in\rho(H)$.

The necessity of~\eqref{estimate1} follows from~\cite{DMS} and the
relation
\[
(R_{\lambda}e_{j,0})_{i,k_i-1}=\begin{cases}
       P_i(\lambda)(Q_j(\lambda)+m(\lambda)P_j(\lambda)),& i\le j;\\
       P_j(\lambda)(Q_i(\lambda)+m(\lambda)P_i(\lambda)),& i>j,\\
 \end{cases}
\]
which directly follows from~\eqref{f_2_5}.
\end{proof}
\begin{remark}
In the case of nonsymmetric tridagonal operators
Theorem~\ref{main_result} was proved in~\cite{AptKV}. 
In fact, we have extended the scheme proposed in~~\cite{AptKV} to 
the case of generalized Jacobi matrices.
A similar result for banded matrices with nonvanishing
extreme diagonals in terms of the corresponding vector polynomials
was obtained in~\cite{BO}.
\end{remark}

\section{The Floquet theory}

In the present section, by using Theorem~\ref{main_result}, we
give a description of spectra of periodic generalized Jacobi
operators.
\begin{definition}\label{periodic_gjm}
Let $s\in\dN$. A generalized Jacobi matrix satisfying the
properties
\[
A_{js+k}=A_k,\quad B_{js+k}=B_{k},\quad\e_{js+k}=\e_k,\quad
j\in\dZ_+,\quad k\in\{0,\dots,s-1\}
\]
will be called an $s$-periodic generalized Jacobi matrix. The
corresponding generalized Jacobi operator in $\cH_{[0,\infty)}$
will be also called an $s$-periodic generalized Jacobi operator.
\end{definition}

Evidently, any $s$-periodic generalized Jacobi matrix satisfies
the assumptions~\eqref{assumption1},~\eqref{assumption2} and we
have
\begin{equation}\label{per_mat}
\cW_{js+k}(\lambda)=\cW_k(\lambda),\quad j\in\dZ_+,\quad
k\in\{0,\dots,s-1\}.
\end{equation}
The main tool for analysis of periodic generalized Jacobi
operators is the following matrix
\[
T(\lambda):=\cW_{[0,s-1]}(\lambda)=\begin{pmatrix}
  -\e_{s-1}b_{s-1} Q_{s-1}(\lambda) & -Q_{s}(\lambda) \\
  \e_{s-1}b_{s-1} P_{s-1}(\lambda) & P_{s}(\lambda) \\
\end{pmatrix}.
\]
The matrix $T(\lambda)$ is called {\it the monodromy matrix}.
Using~\eqref{per_mat}, we get the following relation
\begin{equation}\label{f_3_2}
\cW_{[0,js+k-1]}(\lambda)=T^j(\lambda)\cW_{[0,k-1]}(\lambda),
\quad j\in\dZ_+,\quad k\in\{1,\dots,s\}.
\end{equation}
Let $w_1=w_1(\lambda)$ and $w_2=w_2(\lambda)$ be the roots of the
characteristic equation $\det(T(\lambda)-w)=0$. Introduce the
following notations
\[
E:=\{\lambda\in\dC: |w_1(\lambda)|=|w_2(\lambda)|\},\quad
E_p:=\{\lambda\in\dC:
P_{s-1}(\lambda)=0,\,|b_{s-1}Q_{s-1}(\lambda)|>|P_{s}(\lambda)|\}.
\]
Now, we are ready to give a description of spectra of periodic
generalized Jacobi operators.
\begin{theorem}\label{specper}
The spectrum of an $s$-periodic generalized Jacobi operator has
the form
\[
\sigma(H)=E\cup E_p,\quad \sigma_p(H)=E_p.
\]
\end{theorem}

\begin{proof}
Since $\det T(\lambda)\equiv 1$, we have that
\[
w_1(\lambda)w_2(\lambda)=1.
\]

\noindent{\it Step 1.} First, let us prove that
\[
\{\lambda\in\dC: P_{s-1}(\lambda)\ne
0,\,|w_1(\lambda)|\ne|w_2(\lambda)|\}\subset\rho(H).
\]
To be definite, assume that $|w_1(\lambda)|>|w_2(\lambda)|$. In
this case, we see that
\begin{equation}\label{f_3_2s}
T(\lambda)=
\begin{pmatrix}
x_1&x_2\\
x_3&x_4\\
\end{pmatrix}
\begin{pmatrix}
w_1(\lambda)&0\\
0&w_2(\lambda)\\
\end{pmatrix}
\begin{pmatrix}
x_4&-x_2\\
-x_3&x_1\\
\end{pmatrix},\quad
\det\begin{pmatrix}
x_1&x_2\\
x_3&x_4\\
\end{pmatrix}=1.
\end{equation}
Now, \eqref{f_3_2} can be rewritten in the form
\begin{equation}\label{f_3_3}
\cW_{[0,js+k-1]}(\lambda)=
\begin{pmatrix}
x_1&x_2\\
x_3&x_4\\
\end{pmatrix}
\begin{pmatrix}
w_1^j&0\\
0&w_2^j\\
\end{pmatrix}
\begin{pmatrix}
x_4&-x_2\\
-x_3&x_1\\
\end{pmatrix}
\cW_{[0,k-1]}(\lambda).
\end{equation}
Further, \eqref{f_3_3} is reduced to the following relation
\begin{equation}\label{f_3_4}
\cW_{[0,js+k-1]}(\lambda)=
\begin{pmatrix}
x_1x_4w_1^j-x_2x_3w_2^j & -x_1x_2w_1^j+x_1x_2w_2^j\\
x_3x_4w_1^j-x_4x_3w_2^j & -x_3x_2w_1^j+x_1x_4w_2^j\\
\end{pmatrix}
\cW_{[0,k-1]}(\lambda).
\end{equation}
Multiplying~\eqref{f_3_4} by the vector $(1\,0)^{\top}$ we obtain
\[
\left(%
\begin{array}{c}
  -\e_{k-1}b_{k-1} Q_{js+k-1}(\lambda) \\
  \e_{k-1}b_{k-1} P_{js+k-1}(\lambda) \\
\end{array}%
\right)=
\begin{pmatrix}
x_1x_4w_1^j-x_2x_3w_2^j & -x_1x_2w_1^j+x_1x_2w_2^j\\
x_3x_4w_1^j-x_4x_3w_2^j & -x_3x_2w_1^j+x_1x_4w_2^j\\
\end{pmatrix}
\left(%
\begin{array}{c}
  -\e_{k-1}b_{k-1} Q_{k-1}(\lambda) \\
  \e_{k-1}b_{k-1} P_{k-1}(\lambda) \\
\end{array}
\right).
\]
Thus, for the polynomials $Q_{\cdot}$ one has
\begin{equation}\label{f_3_5}
Q_{js+k-1}(\lambda)=-w_1^j(x_1x_4Q_{k-1}(\lambda)+x_1x_2P_{k-1}(\lambda))+w_2^j(x_2x_3Q_{k-1}(\lambda)+x_1x_2P_{k-1}(\lambda)).
\end{equation}
Similarly, for the polynomials $P_{\cdot}$ we have
\begin{equation}\label{f_3_6}
P_{js+k-1}(\lambda)=-w_1^j(x_3x_4Q_{k-1}(\lambda)+x_3x_2P_{k-1}(\lambda))+w_2^j(x_4x_3Q_{k-1}(\lambda)+x_1x_4P_{k-1}(\lambda)).
\end{equation}
Note that $x_3\ne 0$. Indeed,  if $x_3=0$ then according
to~\eqref{f_3_2s} we would have that $P_{s-1}(\lambda)=0$. Now,
formulas~\eqref{f_3_5} and~\eqref{f_3_6} yield
\begin{equation}\label{f_3_7}
Q_{js+k-1}(\lambda)-\frac{x_1}{x_3}P_{js+k-1}=-w_2^j\left(\frac{x_1}{x_3}P_{k-1}(\lambda)+Q_{k-1}(\lambda)\right).
\end{equation}
For brevity, define $C_k(\lambda):=-\left(\displaystyle{\frac{x_1}{x_3}}P_{k-1}(\lambda)+Q_{k-1}(\lambda)\right)$. 
Then the equality~\eqref{f_3_7} can be rewritten in the following way
\begin{equation}\label{f_3_8s}
Q_{js+k-1}(\lambda)-\frac{x_1}{x_3}P_{js+k-1}=(\widetilde{w}_2)^{js+k-1}C_k(\lambda),\quad
\widetilde{w}_2=w_2^{1/s},\quad k\in\{1,\dots,s\}.
\end{equation}
Since $|w_2|<1$, it follows from~\eqref{f_3_8s} that
$\displaystyle{\xi(\lambda)-\frac{x_1}{x_3}}\pi(\lambda)\in\ell^2_{[0,\infty)}$.
A linear independence of $\pi(\lambda)$ and $\xi(\lambda)$,
and~\eqref{f_3_6},~\eqref{f_3_5}  imply
$\pi(\lambda)\notin\ell^2_{[0,\infty)}$. Now we are ready to
verify the condition~\eqref{estimate1} of
Theorem~\ref{main_result}. Let us assume that $i<j$. Then
\[
\left|P_i(\lambda)\left[Q_j(\lambda)-\frac{x_1}{x_3}P_j(\lambda)\right]\right|=|(-\widetilde{w}_1^if_1^{(k_1)}(\lambda)+\widetilde{w}_2^if_2^{(k_1)}(\lambda))
\widetilde{w}_2^jC_{k_2}(\lambda)|,
\,\widetilde{w}_1=\widetilde{w}_2^{-1},\, k_1,k_2\in\{1,\dots,s\},
\]
where $k_1\equiv i(\text{mod} s)$ and $k_2\equiv j(\text{mod} s)$.
Since
$(n_j-n_i)/N=\left(\sum\limits_{l=i+1}^{j}k_l\right)/N\le\left(\sum\limits_{l=i+1}^{j}N\right)/N=j-i$
we have that $(n_j-n_i)/N\le j-i$. Thus, one obtains
\[
\left|P_i(\lambda)\left[Q_j(\lambda)-\frac{x_1}{x_3}P_j(\lambda)\right]\right|\le|\widetilde{w}_2|^{j-i}C(\lambda)\le
q^{n_j-n_i}C(\lambda),
\]
where
$C(\lambda)=\sup\limits_{k_,k_2\in\{1,\dots,s\}}\{(|f_1^{(k_1)}(\lambda)|+|f_2^{(k_2)}(\lambda)|)|C_{k_2}(\lambda)|\}$
and $q=\widetilde{w}_2<1$.

\noindent {\it Step 2.} Let us show that
\[
\{\lambda\in\dC:
P_{s-1}(\lambda)=0,\,|b_{s-1}Q_{s-1}(\lambda)|>|P_{s}(\lambda)|\}\subset\sigma_p(H).
\]
In this case the monodromy matrix can be represented as follows
\[
T(\lambda)=\begin{pmatrix}
  -\e_{s-1}b_{s-1} Q_{s-1}(\lambda) & -Q_{s}(\lambda) \\
    0                          & P_{s}(\lambda) \\
\end{pmatrix}=
\begin{pmatrix}
1&x_1\\
0&1\\
\end{pmatrix}
\begin{pmatrix}
w_1(\lambda)&0\\
0&w_2(\lambda)\\
\end{pmatrix}
\begin{pmatrix}
1&-x_1\\
0&1\\
\end{pmatrix},
\]
where $w_1(\lambda)=-\e_{s-1}b_{s-1} Q_{s-1}(\lambda)$ and
$w_2(\lambda)=P_{s}(\lambda)$. Further, we have that
\[
\left(%
\begin{array}{c}
  -\e_{k-1}b_{k-1} Q_{js+k-1}(\lambda) \\
  \e_{k-1}b_{k-1} P_{js+k-1}(\lambda) \\
\end{array}%
\right)=
\begin{pmatrix}
w_1^j & w_1^jx_1-w_2^jx_1\\
0     & w_2^j\\
\end{pmatrix}
\left(%
\begin{array}{c}
  -\e_{k-1}b_{k-1} Q_{k-1}(\lambda) \\
  \e_{k-1}b_{k-1} P_{k-1}(\lambda) \\
\end{array}
\right),
\]
that is, the following formulas hold true
\[
P_{js+k-1}(\lambda)=w_2^jP_{k-1}(\lambda),
\]
\[
Q_{js+k-1}(\lambda)=-w_1^jQ_{k-1}(\lambda)+(w_1^j-w_2^j)x_1P_{k-1}(\lambda).
\]
Since $|w_2(\lambda)|<|w_1(\lambda)|$,
$G\pi(\lambda)\in\ell^2_{[0,\infty)}$ and, therefore,
$\pi(\lambda)\in\ell^2_{[0,\infty)}$. So, we have proved that
$E_p\subset\sigma(H)$.

\noindent{\it Step 3.} By the same reasoning as in the Step 1 it
can be shown that
\[
\{\lambda\in\dC:
P_{s-1}(\lambda)=0,\,|b_{s-1}Q_{s-1}(\lambda)|>|P_{s}(\lambda)|\}
\subset\rho(H).
\]

\noindent{\it Step 4.} We complete this proof by proving the
following inclusion
\[
\{\lambda\in\dC: |w_1(\lambda)|=|w_2(\lambda)|\}\subset\sigma(H).
\]
Since $w_1w_2=1$ we have that $|w_1(\lambda)|=|w_2(\lambda)|=1$.
If $T$ is diagonalizable  then, due to~\eqref{f_3_5},
\eqref{f_3_6}, there exist numbers $\alpha_k^{(1)}(\lambda),
\alpha_k^{(2)}(\lambda),\beta_k^{(1)}(\lambda),\beta_k^{(2)}(\lambda)\in\dC$
such that
\begin{equation}\label{f_3_8}
Q_{js+k-1}(\lambda)=w_2^{-j}\alpha_k^{(1)}+w_2^j\alpha_k^{(2)},\quad
j\in\dZ_+,\quad k\in\{1,\dots,s\}.
\end{equation}
\begin{equation}\label{f_3_9}
P_{js+k-1}(\lambda)=w_2^{-j}\beta_k^{(1)}+w_2^j\beta_k^{(2)},\quad
j\in\dZ_+,\quad k\in\{1,\dots,s\}.
\end{equation}
Since the sequences $Q_{\cdot}$ and $P_{\cdot}$ are linearly
independent, the sequence $\{Q_j+mP_j\}_{j=0}^{\infty}$ is nonzero
for any $m\in\dC$. Also, observe that the sequence
$\{|w^j\alpha+\beta|\}_{j=1}^{\infty}$ ($\alpha\ne 0$, $|w|=1$)
converges only for $w=1$. Taking into account these observations
one can conclude that the series
\[
\sum_{j=0}^{\infty}|Q_j(\lambda)+mP_j(\lambda)|^2=
\sum_{i=0}^{\infty}\sum_{k=0}^{s-1}|\alpha_k^{(1)}+m\beta_k^{(1)}+w_2^{2i}(\alpha_k^{(2)}+m\beta_k^{(2)})|^2
\]
does not converge for any number $m\in\dC$. So, we have that
$\xi(\lambda)+m\pi(\lambda)\notin\ell^2_{[0,\infty)}$ for any
$m\in\dC$. Similarly, one can conclude that
$\pi(\lambda)\notin\ell^2_{[0,\infty)}$.

If $T$ is similar to a Jordan block then $w_1=w_2$ and $w_1w_2=1$.
So, $w_1=\pm 1$. First, let us consider the case $w_1=1$. The
monodromy matrix takes the form
\[
T(\lambda)=
\begin{pmatrix}
x_1&x_2\\
x_3&x_4\\
\end{pmatrix}
\begin{pmatrix}
1&1\\
0&1\\
\end{pmatrix}
\begin{pmatrix}
x_4&-x_2\\
-x_3&x_1\\
\end{pmatrix}.
\]
Further,~\eqref{f_3_2} yields
\[
\cW_{[0,js+k-1]}(\lambda)=
\begin{pmatrix}
x_2x_4-x_1x_3j-x_2x_3 & -x_2^2+x_1^2j-x_1x_2\\
-x_3^2j               & -x_3x_2+x_3x_1j+x_4x_1\\
\end{pmatrix}
\cW_{[0,k-1]}(\lambda).
\]
Thus for any $m\in\dC$ the vectors $\xi(\lambda)+m\pi(\lambda)$
and $\pi(\lambda)$ do not belong to $\ell^2_{[0,\infty)}$.
Analogously, one can consider the case $w_1=-1$.
\end{proof}
\begin{remark}
A description of spectra of classical Jacobi operators was
obtained in~\cite{Ger}.
\end{remark}
Let us remind that $w_1$ and $w_2$ are the roots of the equation
\begin{equation}\label{f_3_10}
w^2-(P_s(\lambda)-\e_{s-1}b_sQ_{s-1}(\lambda))w+1=0.
\end{equation}
Due to~\eqref{f_3_10}, one has
\begin{equation}\label{f_3_11}
w_1+w_2=P_s(\lambda)-\e_{s-1}b_sQ_{s-1}(\lambda),
\end{equation}
\begin{equation}\label{f_3_12}
w_1w_2=1.
\end{equation}
\begin{remark}[\cite{BK}]\label{rem_3_3}
Formulas~\eqref{f_3_11},~\eqref{f_3_12} allow us to give another
description of $E$
\begin{equation}\label{h_ex_1}
E=\{\lambda\in\dC:
(P_s(\lambda)-\e_{s-1}b_sQ_{s-1}(\lambda))\in[-2,2]\}.
\end{equation}
\end{remark}

\begin{example} Let us consider the following difference equation
\begin{equation}\label{ex_1}
\frac{1}{2}u_{j+1}-\lambda^2u_j-\frac{1}{2}u_{j-1}=0\quad
(j\in\dN).
\end{equation}
Clearly, the three-term recurrence relations~\eqref{ex_1} generate
the following $1$-periodic generalized Jacobi matrix
\[
H=\begin{pmatrix}
A   & B & {\bf 0}\\
B   & A &\ddots    \\
{\bf 0} &  \ddots&\ddots\\
\end{pmatrix},\quad
A=\begin{pmatrix}
0&0\\
1&0\\
\end{pmatrix},\quad
B=\begin{pmatrix}
0&\frac{1}{2}\\
0&0\\
\end{pmatrix}.
\]
In this case, by easy calculations, we have
\[
P_0(\lambda)=1,\quad P_1(\lambda)=2\lambda^2,\quad
Q_0(\lambda)=0,\quad Q_1(\lambda)=2.
\]
Since $P_0(\lambda)=1\ne 0$, it follows from Theorem~\ref{specper}
that $\sigma_p(H)=\emptyset$. Next, according to
Theorem~\ref{specper} and~\eqref{h_ex_1}, we have
\[
\sigma(H)=\{\lambda\in\dC: 2\lambda^2\in[-2,2]\}=[-1,1]\cup[-i,i].
\]
\end{example}

 By the same reasoning as in~\cite{BK}, one can prove the
following statement on the structure of the set $E$.
\begin{proposition}[\cite{BK}]\label{structure}
The compact set $E$ has no interior points. The open set
$D:=\dC\setminus E$ is connected. The functions $w_1$ and $w_2$
are single-valued in  $D$.
\end{proposition}

\section{Pad\'e approximants}

Our goal in this section is to prove convergence results for
Pad\'e approximants.

\begin{definition}[\label{Pade}{\rm\cite{NikSor}}] The $[L/M]$ Pad\'e approximant
to the function
$\f(\lambda)=\displaystyle{-\sum\limits_{j=0}^{+\infty}\frac{s_j}{{\lambda}^{j+1}}}$
is defined as a ratio
\[
f^{[L/M]}(\lambda)=\frac{A^{[L/M]}\left(\frac{1}{\lambda}\right)}
{B^{[L/M]}\left(\frac{1}{\lambda}\right)}
\]
of two polynomials $A^{[L/M]}$, $B^{[L/M]}$ of formal degree $L$
and $M$, respectively, such that $B^{[L/M]}(0)\ne 0$ and
\[
\sum_{j=0}^{+\infty}\frac{s_j}{{\lambda}^{j+1}}+f^{[L/M]}(\lambda)=
O({\lambda}^{-(L+M+1)})\quad(\lambda\to\infty).
\]
\end{definition}

In the case $L=M=n$, the $[n/n]$ Pad\'e approximant is also called
{\it the $n$-th diagonal Pad\'e approximant}.

Let us consider the Weyl function
$m(\lambda)=[(H-\lambda)^{-1}e,e]$, where the corresponding matrix
$H$ satisfies the assumptions~\eqref{assumption1},
\eqref{assumption2}. The representation~\eqref{SWfun} of
$m_{[0,j-1]}$ yields
\begin{equation}\label{Qasympt-1}
m_{[0,j-1]}(\lambda)=- \frac{Q_j(\lambda)}{P_j(\lambda)}=
-\frac{s_{0}}{\lambda}-\frac{s_{1}}{\lambda^{2}}-\dots-\frac{s_{2n_j-2}}{\lambda^{2n_j-1}}
+O\left(\frac{1}{\lambda^{2n_j}}\right) \quad(\lambda\to\infty),
\end{equation}
where ${s_{i}}=[H^ie,e]$ ($i=0,\dots,2n_j-2$). Moreover, it was
shown in~\cite{DD} that  $m_{[0,j-1]}(\lambda)$ has the following
asymptotic expansion
\begin{equation}\label{Qasympt+1}
m_{[0,j-1]}(\lambda)=
-\sum_{i=0}^{2n_j-2+k_j}\frac{s_{i}}{\lambda^{i+1}}
+O\left(\frac{1}{\lambda^{2n_j+k_j}}\right)
\quad(\lambda\to\infty),
\end{equation}
where ${s_{i}}=[H^ie,e]$ ($i=0,\dots,2n_j-2+k_j$). The latter
means that the rational function
\begin{equation}\label{PadeAp}
f^{[n_j/n_j]}(\lambda)=m_{[0,j-1]}(\lambda)=\frac{\displaystyle
A^{[n_j/n_j]}\left({1}/{\lambda}\right)} {\displaystyle
B^{[n_j/n_j]}\left({1}/{\lambda}\right)}=
\frac{\displaystyle-\frac{1}{\lambda^{n_j}}Q_j(\lambda)}{\displaystyle\frac{1}{\lambda^{n_j}}P_j(\lambda)}\quad
(j=1,2,\dots)
\end{equation}
is the $[n_j/n_j]$ Pad\'e approximant to  $m$. Besides, it follows
from the Pad\'e theorem (see~\cite[Theorem~1.4.3]{Baker}), that
for $L$ and $M$ satisfying
\[
L\ge n_j,\,\,\,M\ge n_j,\,\,\, L+M\le 2n_j+{k_j}-1
\]
the $[L/M]$ Pad\'e approximants coincide with $f^{[n_j/n_j]}$, and
for $L$, $M$ satisfying
\[
L\le n_j+{k_j}-1,\,\,\,M\le n_j+{k_j}-1,\,\,\, L+M\ge 2n_j+{k_j}
\]
the $[L/M]$ Pad\'e approximants do not exist  (for details
see~\cite{DD}).

In what follows we need the following definition.

\begin{definition}\label{num_ran}
The set $\Theta(H):=\{(Hy,y)_{\ell^2}: \Vert y\Vert=1\}\subset\dC$
is called a numerical range of the operator $H$.
\end{definition}
Clearly, the numerical range of a bounded operator is a bounded
set. By the Hausdorff theorem we have that
$\sigma(H)\subset\overline{\Theta(H)}$ (see~\cite{Kat}).
\begin{theorem}\label{un_loc_con}
Let $m(\lambda)=[(H-\lambda)^{-1}e,e]$ be the Weyl function of $H$
satisfying~\eqref{assumption1}, \eqref{assumption2}. Then there
exists a subsequence of diagonal Pad\'e approximants
$f^{[n_j/n_j]}$ to $m$, which converges to $m$ locally uniformly
in $\dC\setminus\overline{\Theta(H)}$.
\end{theorem}
\begin{proof}
First, note that $\Theta(H_{[0,n]})\subset\Theta(H)$, and,
therefore,
$\overline{\Theta(H_{[0,n]})}\subset\overline{\Theta(H)}$. As a
consequence, we have that if
$\lambda\in\dC\setminus\overline{\Theta(H)}$ then
$\lambda\in\rho(H_{[0,n]})$ for all $n\in\dZ_+$. Let
$\lambda\in\dC\setminus\overline{\Theta(H)}$. Then
$\ran(H-\lambda)=\ell^2_{[0,\infty)}$ and for any finite vector
$\phi$ we have
\[
(H_{[0,j]}-\lambda)^{-1}\phi\to(H-\lambda)^{-1}\phi,\quad
j\to+\infty.
\]
Thus, one can obtain
\begin{equation}\label{conv_pa}\begin{split}
f^{[n_j/n_j]}(\lambda)=m_{[0,j-1]}(\lambda)=[(H_{[0,j-1]}-\lambda)^{-1}e,e]\\=((H_{[0,j-1]}-\lambda)^{-1}e,Ge)_{\ell^2}\to((H-\lambda)^{-1}e,Ge)_{\ell^2}=m(\lambda).\end{split}
\end{equation}
According to~\cite[p.336 Theorem 3.2]{Kat} one has
\[
\Vert(H_{[0,j]}-\lambda)^{-1}\Vert\le
\frac{1}{\dist(\lambda,\Theta(H_{[0,n]}))}
\le\frac{1}{\dist(\lambda,\Theta(H))}.
\]
Hence, the family  $f^{[n_j/n_j]}=m_{[0,j]}$ is uniformly bounded
on compact sets in $\dC\setminus\overline{\Theta(H)}$. Actually,
the following estimates hold true
\[
|m_{[0,j]}(\lambda)|=|((H_{[0,j]}-\lambda)^{-1}e,Ge)_{\ell^2}|\le\Vert(H_{[0,j]}-\lambda)^{-1}e\Vert_{\ell^2}\Vert
Ge\Vert_{\ell^2}\le \frac{\Vert
Ge\Vert_{\ell^2}}{\dist(\lambda,\Theta(H))}.
\]
So, the family  $f^{[n_j/n_j]}$ is uniformly bounded and,
therefore, $f^{[n_j/n_j]}$ is precompact. Now, to complete the
proof it is sufficient to apply~\eqref{conv_pa} and the Vitali
theorem.
\end{proof}
\begin{remark}
In the proof of Theorem~\ref{un_loc_con} we used the method
proposed in~\cite{BK} for complex Jacobi matrices. Note that
Theorem~\ref{un_loc_con} is a generalization of~\cite[Theorem
1]{Gonc82}. More precisely, we do not suppose existence of all
diagonal Pad\'e approximants and all poles of the existed diagonal
Pad\'e approximants belong to the convex set
$\overline{\Theta(H)}$.
\end{remark}

 Further, following the scheme proposed in~\cite{AptKV}, we find
out a behavior of  the associated polynomials at the points of the
resolvent set.

\begin{proposition}\label{pol_as}
For any $\lambda\in\rho(H)$ the following inequality holds true
\begin{equation}\label{n_pol_as}
\limsup_{j\to+\infty}|P_j(\lambda)|^{1/j}>1.
\end{equation}
\end{proposition}
\begin{proof} In fact, the proof is in line with that in~\cite{AptKV}.
However, we give the proof here for the convenience of the
readers. Using~\eqref{Ostrogr2} and~\eqref{reprWs}, one can see
\begin{equation}\label{PA_h_1}
P_j(\lambda)(\e_j b_j W_{j+1}(\lambda))-P_{j+1}(\lambda)(\e_j b_j
W_j(\lambda))=1.
\end{equation}
Due to~\eqref{assumption1}, \eqref{assumption2},
and~\eqref{estimate1}, one has
\begin{equation}\label{PA_h_2}
|\e_j b_j W_{j}(\lambda)|\le C_1{q_1}^{j},\quad 0<q_1<1.
\end{equation}
It follows from~\eqref{PA_h_1} and~\eqref{PA_h_2} that the
sequence $P_j(\lambda)$ can not be majorized by a geometric
sequence $p^j$ with $p<1/q_1$, that is, for any $p<1/q_1$ and any
positive constant $C_2$ the inequality
\[
|P_j(\lambda)|\le C_2p^j
\]
is not satisfied for an infinite number of indices $j$. We thus
have
\[
|P_{j_k}(\lambda)|\ge C_2p^{j_k}.
\]
If we choose $p$ such that $1/q_1>p>1$ we obtain the required
result.
\end{proof}

From these statements one can deduce the following result on
convergence of diagonal Pad\'e approximants.

\begin{theorem}\label{pointwise}
Under the assumptions~\eqref{assumption1}, \eqref{assumption2},
for any $\lambda\in\rho(H)$ there exists a subsequence of diagonal
Pad\'e approximants to  $m(\lambda)=[(H-\lambda)^{-1}e,e]$, which
converges to $m(\lambda)$ at $\lambda$.
\end{theorem}
\begin{proof}
From~\eqref{n_pol_as} we get that there exists a subsequence $j_k$
such that
\begin{equation}\label{f_subs}
|P_{j_k}(\lambda)|\ge C_2p^{j_k},\quad p>1.
\end{equation}
Further, Theorem~\ref{main_result} and~\eqref{f_subs} imply the
relation
\[
\left|m(\lambda)+\frac{Q_{j_k}(\lambda)}{P_{j_k}(\lambda)}\right|=
\left|\frac{W_{j_k}(\lambda)}{P_{j_k}(\lambda)}\right|\le
C_3\left(\frac{q}{p}\right)^{j_k},
\]
which completes the proof.
\end{proof}

\noindent{\bf Acknowledgments}. I express my gratitude to
Professor V.A.~Derkach for reading the manuscript and giving
helpful comments. I also thank the anonymous referee for helpful remarks.

\end{document}